\documentclass[a4paper,twoside,11pt,reqno]{amsart}

\usepackage[scaled=0.92]{helvet}
\usepackage{courier}
\usepackage[T1]{fontenc}
\usepackage[english]{babel}
\usepackage{amsthm,amsfonts,amssymb}
\usepackage{hyperref}
\usepackage{mathrsfs}

\hypersetup{
colorlinks=true,
citecolor=red,
linkbordercolor={1 1 1},
citebordercolor={1 1 1},
}
\theoremstyle{plain}
\newtheorem{satz}{Theorem}[section]
\newtheorem{prop}[satz]{Proposition}
\newtheorem{cor}[satz]{Corollary}
\newtheorem{lem}[satz]{Lemma}

\theoremstyle{definition}
\newtheorem{defn}[satz]{Definition}
\newtheorem{rem}[satz]{Remark}
\newtheorem{cond}[satz]{Condition}

\newcommand{\rw}{\rightarrow}
\newcommand{\de}{\displaystyle}
\newcommand{\ep}{\varepsilon}
\newcommand{\ml}{\mathcal}
\newcommand{\tf}{\textbf}
\newcommand{\eh}{\emph}
\newcommand{\im}{\item}
\newcommand{\ue}{\underline}

\newcommand{\pl}{\partial}
\newcommand{\ph}{\varphi}
\newcommand{\nt}{\noindent}

\newcommand{\x}{\times}

\newcommand{\R}{\mathbb{R}}

\newcommand{\Z}{\mathbb{T}}
\newcommand{\N}{\mathbb{N}}
\newcommand{\T}{\mathbb{T}}

\newcommand{\norm}[1]{\lVert#1\rVert}

\DeclareMathOperator{\rk}{rk}

\DeclareMathOperator{\diam}{diam}
\DeclareMathOperator{\diag}{diag}
\DeclareMathOperator{\modulus}{mod}

\begin{document}

\title{Fast drift and diffusion in an example of isochronous system through windows method}
\author{Alessandro Fortunati} 
\email{alessandro.fortunati2012@gmail.com}
	
\date{}
\begin{abstract}
We study the problem of Arnold's diffusion in an example of isochronous system by using a geometrical method known as Windows Method. Despite the simple features of this example, we show that the absence of an anisochrony term leads to several substantial difficulties in the application of the method, requiring some additional devices as non-equally spaced transition chains and variable windows. In this way we are able to obtain a set of fast orbits whose drifting time matches, up to a constant, the time obtained via variational methods.  
\end{abstract}

\maketitle

\section{Introduction}\label{sec:intro}
Arnold's Diffusion \cite{arn} is a topological instability phenomenon arising in nearly integrable Hamiltonian systems with more than two degrees of freedom. It consists on the existence of a class of motions whose effect is to produce a drift of slow variables of order one in a ``very long'' time $T_d$, i.e. tending to infinity as the perturbation size $\mu$ approaches to zero, for \ue{all} sufficiently small $\mu$. A quite relevant problem clearly concerns the estimate of the minimal time in order to observe this phenomenon. The earlier literature on this topic, such as \cite{chigal94}, \cite{bes} etc., has shown that the attempts to give an estimate of $T_d$ crucially depend on the method used to construct such class of trajectories, rather than the features of the system at hand, as it would be natural to expect. The question is addressed in \cite{gal97} and in \cite{gal99}, by taking as a paradigmatic example the following simple a priori-unstable isochronous system 
\begin{equation}\label{eq:hamiltonian}
H(p,q,\ue{A},\ue{\ph})=\ue{\omega} \cdot
\ue{A}+\frac{p^2}{2}+(\cos q-1)+\mu F(q,\ue{\ph}) \mbox{,}
\end{equation}
with $(q,\ue{\ph},p,\ue{A}) \in \ml{M}:=T^* (\mathbb{S}^1 \times \mathbb{T}^2)$. Some relevant subsequent works on this class of systems (with some various generalizations) 
extent this gap further, as summarized in the following (non-exhaustive) table\\
\begin{table}[h!]
\begin{center}
\begin{tabular}{|l|l|l|}
\hline
Work &  $T_d$ & approach\\
\hline
\hline
\cite{gal97} & super-exponential & geometrical \\
\hline
\cite{gal99} & exponential (i.e. $O(\exp(\mu^{-1}))$)& geometrical \\
\hline
\cite{crescras} & polynomial (i.e. $O(\mu^{-\gamma})$)& geometrical \\
\hline
\cite{berbol} & logarithmic (i.e. $O(\mu^{-1}\log \mu^{-1})$)& variational \\
\hline
\end{tabular}
\end{center}
\end{table}
\\
where $\gamma$ is a positive $O(1)$ constant.\\In the much more investigated anisochronous case, the situation was a little bit more encouraging. In particular, the logarithmic estimate finally given by the variational approach of \cite{berbiasbol}, was also obtained in \cite{cresguil3} via the so-called Windows Method (WM). This powerful geometrical tool (originally due to Aleskeev and Easton, see for instance \cite{east81} and references therein), was reconsidered and developed in \cite{mar96} on a simplified version of the Arnold's example. Nowadays, the WM has been deeply improved leading to some advanced applications on the problem of Arnold's diffusion, see for instance, \cite{gidzgli}, \cite{gidrob} and \cite{lara}.\\ 
The motivation of this work arises from the fact that a treatment of the isochronous case cannot be obtained neither from \cite{mar96} nor from \cite{cresguil3} by taking as zero the anisochrony coefficient. The main reason can be briefly described as follows. As done in \cite{cres3}, it is possible to show the existence of a conjugacy between the dynamics in a neighbourhood of a Graff torus and a symbolic dynamics on a suitable alphabet and then, in particular, the existence of periodic orbits. This allows to construct in \cite{cresguil3}, a chain of (hyperbolic) periodic orbits $\ml{O}_k$, surrounding a transition chain of partially hyperbolic tori. Arnold's diffusion is obtained via WM as a shadowing of $\ml{O}_k$.\\
It is clear that the above described argument cannot be applied in the isochronous case as this class of systems does not admit periodic orbits: the (Diophantine) frequencies of rotators are fixed and then trajectories are open for all $t$, densely filling the underlying flat torus. Beyond the
impossibility in adapting the construction done in \cite{cresguil3}, a deeper obstruction arises in the isochronous case: the absence of an anisochrony term inhibits a remarkable geometrical phenomenon pointed out in  \cite{cresguil1} and called ``transversality-torsion''. The latter is related to a suitable compressing-stretching action of the phase flow on windows, giving a key condition for a shadowing result (called \emph{correct alignment}) in a quite natural way. This cannot be done in our case and the correct alignment condition requires a more careful treatment.\\
In order to overcome this difficulty, we have considered a suitable sequence of ``self deforming'' (i.e. non-constant) windows, calling this mechanism \emph{simulated torsion}. Such a sequence has been determined by solving a one dimensional discrete dynamical system over the transition chain.\\The more restrictive estimates (of technical nature) arising in the isochronous case, also obstruct the use of an equally spaced transition chain\footnote{namely, a transition chain in which the distance in the actions space of two consecutive tori is constant.}. For this purpose, we have used the \emph{elastic chain} tool, introduced for the first time in \cite{gal99} and consisting in a step-by-step variation of an equally spaced chain\footnote{The use of elastic chains is a key ingredient for a remarkable reduction of the transition time in the geometrical methods context. A review on the differences between geometrical approaches of \cite{chigal94} and of \cite{east81} and how these are able to improve the speed of diffusion, can be found in \cite{for}.}. 
In this way, we are able to suitably move the hetero/homoclinic points, in order to bypass the ergodization time\footnote{i.e. the time necessary to approach a given point on the flat torus within a prefixed distance.}. This leads to a systematic reduction of the transition time allowing us to obtain a logarithmic estimate for $T_d$.

\section{Preliminary facts and main result}\label{sec:syst}
\subsection{System at hand, invariant tori and whiskers intersection}
Let us consider the system (\ref{eq:hamiltonian}) where the perturbing function is supposed of the form $F(q,\ue{\ph})=(1-\cos q) f(\ue{\ph})$, with
\begin{equation}\label{eq:perturbazione}
f(\underline{\ph})=\sum_{\ue{k} \in \mathbb{Z}^2  \atop |k_1|+|k_2|
\leq \Gamma} f_{\ue{k}} \cos (\ue{k} \cdot \ue{\ph}) \mbox{,}
\end{equation}
$\Gamma$ is a fixed positive constant, and $f_{\ue{k}}=f_{-\ue{k}}>0$ for all
$\ue{k}$. In particular $f$ is an even trigonometric polynomial. We suppose $\ue{\omega}$ to be a $(C,\tau)-$Diophantine vector\footnote{There exist constants $C \neq 0$ and $\tau \geq 1$ such that $|\ue{\omega}\cdot \ue{k}| \geq
C|\ue{k}|^{-\tau}$.}.\\
\noindent First of all note that the unperturbed system possesses a 
two parameter continuous family of invariant tori 
\[ \ml{T}(\ue{A})=\{(p,q,\ue{A},\ue{\ph}):p=q=0,\quad \ue{\ph} \in \T^2\} \mbox{,} \]
with $\ue{A} \in \R^2$. Due to term $1-\cos q$ in the perturbing function, the entire family survives also in the perturbed system. This is of course one of the main advantages of this simple example. If one requires that two tori with different actions $\ue{B},\ue{C}$ lie on the same energy level, these have to satisfy  $\ue{\omega} \cdot (\ue{B}-\ue{C})=0$. Hence we have
\begin{equation}\label{eq:isoenergy}
\ue{B}=\ue{C}+\delta \ue{\omega}^{\perp} \mbox{,}
\end{equation}
for some $\delta$. So they form a one parameter continuous family.\\ 
These tori are \emph{whiskered} in the sense that each of them admits a couple of invariant manifolds $W^{s,u}(\ml{T}(\ue{A}))$ called \emph{whiskers}, such that, by denoting with $\Phi^t$ the phase flow given by (\ref{eq:hamiltonian}), $ \lim_{t \rw \infty} \Phi^t(\ue{z})   \in \ml{T}(\ue {A}) $ for all $\ue{z} \in W^s (\ml{T}(\ue{A}))$ and $ \lim_{t \rw -\infty} \Phi^t(\ue{z})   \in \ml{T}(\ue {A}) $ for all $\ue{z} \in W^u (\ml{T}(\ue{A}))$.\\
As can be easily seen by equations of motion, the whiskers (positive branches) $W^u(\mathcal{T}(\ue{B}))$ and
$W^s(\mathcal{T}(\ue{C}))$ possess the following structure in the unperturbed system
$$
\begin{array}{l}
W^u(\mathcal{T}(\ue{B})) =\{(q,\ue{\ph},p,\ue{A}):p_0(q):=
\sqrt{2(1-\cos q)},\quad \ue{A}=\ue{B}, \quad (q,\ue{\ph}) \in (0,2
\pi) \x
\mathbb{T}^2\}\vspace{0.1cm}\\
W^s(\mathcal{T}(\ue{C})) =\{(q,\ue{\ph},p,\ue{A}):p_0(q):=
\sqrt{2(1-\cos q)},\quad \ue{A}=\ue{C}, \quad (q,\ue{\ph}) \in (0,2
\pi) \x \mathbb{T}^2\}
\end{array}
$$
It is clear that these  whiskers intersect each other if, and only if, 
$\delta=0$ and moreover the intersection is flat, meaning that the whiskers tangent spaces coincide at every point. Note that the whiskers possess the structure (trivial in the $\ue{A}$) of graphs over the angles.\\
By continuous dependence on parameters, the whiskers will keep the graph structure also in the perturbed case, for sufficiently small $\mu$. So, for all positive $O(1)$ constants $d$, the whiskers can be written as 
\begin{equation}\label{eq:baffiperturbati}
\begin{array}{l}
W^u(\mathcal{T}(\ue{B})) =\{(q,\ue{\ph},p_0(q)+\mu
p^u(q,\ue{\ph},\mu),\quad \ue{B}+\mu \ue{A}^u(q,\ue{\ph},\mu)),
(q,\ue{\ph}) \in (0, 2 \pi-d) \x
\mathbb{T}^2\} \vspace{0.1cm}\\
W^s(\mathcal{T}(\ue{C})) =\{(q,\ue{\ph},p_0(q)+\mu
p^s(q,\ue{\ph},\mu),\quad \ue{C}+\mu \ue{A}^s(q,\ue{\ph},\mu)),
(q,\ue{\ph}) \in (d,2 \pi) \x \mathbb{T}^2\}
\end{array}
\end{equation}
see e.g. \cite{gal0}, also for its recursive construction at all orders.\\
Let us denote with
$$\Delta p:=p^u(q,\ue{\ph},\mu)-p^s(q,\ue{\ph},\mu),\qquad
\Delta \ue{A}:=\ue{A}^s(q,\ue{\ph},\mu)-\ue{A}^u(q,\ue{\ph},\mu)
\mbox{,}$$ and with
\begin{equation}\label{eq:fqfidelta}
\ue{F}(q,\ue{\ph},\delta):=\delta \ue{\omega}^{\perp}+\mu \Delta \ue{A}(q,\ue{\ph},\mu) \mbox{.}
\end{equation}
Now define 
\[ M:=\ml{D}\ue{F} \equiv \frac{\pl (\Delta  p,\Delta \ue{A})}{\pl(q,\ue{\ph})} \mbox{.}\]
It is an easy consequence of the isoenergy constraint between $\ml{T}(\ue{B})$ and $\ml{T}(\ue{C})$ that if there exists a point $(q^*,\ue{\ph}^*) \in (d,2 \pi-d) \x
\T^2$ such that $\ue{F}(q^*,\ue{\ph}^*,\delta)=\ue{0}$ and $\rk M=2$, then\footnote{ $\pitchfork$ denotes the transversal intersection.} $W^u(\ml{T}(\ue{B})) \pitchfork W^s(\ml{T}(\ue{C}))$
at $(q^*,\ue{\ph}^*)$. 
Actually, the condition $\rk M=2$ is maximal as by the conservation of energy one gets 
\begin{equation}
\begin{array}{rcl}\label{eq:isoenqphi}
\pl_q \Delta p&=&\alpha \ue{\omega} \cdot \pl_q \Delta \ue{A}\\
\pl_{\ph_k} \Delta p&=&\alpha \ue{\omega} \cdot \pl_{\ph_k} \Delta \ue{A},\qquad k=1,2 \mbox{.}
\end{array}
\end{equation}
with $\alpha:=-1/(p_0(q^*)+\mu p^s(q^*,\ue{\ph}^*,\mu))$, so $\rk M<3$.\\
In this case we shall speak, as usual, of a transversal homoclinic point if $\delta=0$ (i.e. $\ue{B}=\ue{C}$) and of a transversal heteroclinic point otherwise.\\ 
Another key ingredient is the symplectic character of the whiskers\footnote{see \cite{locmarsau} for a comprehensive treatment of this topic and a different algorithmical approach for the whiskers building at all orders in $\mu$.}. By its use we can state the following  
\begin{prop}
Suppose that, by defining as $M_{\mbox{se}}$ the south-east $2 \times 2$ submatrix of $M$, $\det M_{\mbox{se}}(q^*,\ue{\ph}^*)>0$ holds for some $(q^*,\ue{\ph}^*)$. Then we have 
\begin{equation}\label{eq:deltapneqzero}
\pl_q \Delta p \neq 0; \qquad 
\det 
\left(%
\begin{array}{cc}
  \pl_q \Delta A_1 & \pl_{\ph_1} \Delta A_1 \\
  \pl_q \Delta A_2 & \pl_{\ph_1} \Delta A_2 \\
\end{array}%
\right)(q^*,\ue{\ph}^*) \neq 0
\end{equation}
\end{prop}

\proof (Sketch) Due to their lagrangian nature the whiskers are gradient of two respective (family of) generating functions $S^{u,s}(q,\ue{\ph},\mu)$ so $M$ is nothing but the Hessian matrix of $S^u(q,\ue{\ph},\mu)-S^s(q,\ue{\ph},\mu)$, computed at the intersection point. This implies that $M$ is a symmetric matrix. By solving the system given by conditions (\ref{eq:isoenqphi}) with respect to the entries of $M_{\mbox{se}}$ and imposing that $\det M_{\mbox{se}}>0$ the statement easily follows. 
\endproof
\noindent Now we are going to assume the following typical\footnote{the existence of a transversal homoclinic point is a standard result for this system; see e.g. \cite{galgenmas} for the same result in a more general context. Condition $\det M_{\mbox{se}}>0$ (instead of $\det M_{\mbox{se}} \neq 0$) is a feature of this system and can be straightforwardly checked from the Melnikov integral (for sufficiently small $\mu$).}
\begin{cond}[splitting]\label{cond:splitting}
The system (\ref{eq:hamiltonian}) possesses a transversal homoclinic point at $q=\pi$. In other terms there exists a point $(q,\ue{\ph})=(\pi,\ue{\ph}^0)$ such that $F(\pi,\ue{\ph}^0)=0$ and $\det M_{\mbox{se}}(\pi,\ue{\ph}^0)>0$ with $(M_{\mbox{se}})_{ij}=O(1)$.
\end{cond}
\begin{rem}\label{rem:splitting}
Clearly $(\pi,\ue{\ph}^0)$ is a particular homoclinic point. As the whiskers intersection takes place along an entire orbit, all the points of the form 
$(q,\ue{\ph}^0+\ue{\omega}T_{q \rw \pi})$, where $T_{q \rw \pi}$ is the time to evolve from  $\pi$ to $q$, are homoclinic as well. If $q$ is $O(1)$ bounded away from the origin, $M_{\mbox{se}}$ can be well approximated via Melnikov integral. It can be easily seen that $M_{\mbox{se}}$ does not change (up to $O(\mu^2)$) if different section for $q$ are chosen.
\end{rem}
The section $q=\pi$ is useful in presence of symmetries of the perturbing function as in this case. It is possible to show (see e.g. \cite{galgenmas}) that, by parity,  $(q,\ue{\ph}^0)=(\pi,\ue{0})$ is a transversal homoclinic point at all orders in $\mu$. We denote as $\ue{\ph}^0=:(\ph_1^0,\ph_2^0)$.\\
In the current setting, we state the following
\begin{lem}\label{lemma:splitting}
Let us assume the splitting condition, and define the section
$$\Sigma:=\{\ue{\ph} \in \mathbb{T}^2: \ph_2 =\ph_2^0\} \mbox{.}$$
There exist $\bar{\delta}=O(\mu)$ and two functions $q,\ph_1$ defined on $(-\bar{\delta},\bar{\delta})$, with $q(0)=\pi$ and
$\ph_1(0)=\ph_1^0$, such that, for all $\ue{B}$ and sufficiently small $\mu$

$$W^u(\ml{T}(\ue{B}))\pitchfork W^s(\ml{T}(\ue{B}+\delta \ue{\omega}^{\perp})) \mbox{,}$$
on  $\Sigma$ at the point $(q(\delta),\ph_1(\delta))$. In particular
\begin{equation}\label{eq:variabilitaphiuno}
  \frac{d}{d \delta}\ph_1(\delta) \neq 0,
  \qquad \frac{d}{d \delta} \ph_1(\delta), \frac{d}{d \delta} q(\delta)=O(\mu^{-1}) \mbox{.}
\end{equation}
for all $\delta \in (-\bar{\delta},\bar{\delta})$.
\end{lem}
The just stated version of this classical splitting result is a sort of counterpart of the statement which can be found in \cite{gal99}. According to the latter, as $\delta$ varies, the $\ue{\ph}$ coordinates of the intersection point move along a curve on $\T^2$, while $q$ is fixed to $\pi$. In our case $q$ moves, but $\ph_2$ is fix, allowing us to use a single section on $\T^2$ for all values of $\delta$ and then to perform an advantageous dimensional reduction.
\proof
By hypothesis $\ue{F}(\pi,\ue{\ph}^0,0)=\ue{0}$ and $\rk M_{se}=2$, so by implicit function theorem, there exist $\bar{\delta}_1$ and two functions 
$q(\delta),\ph_1(\delta)$ such that $(q(0),\ph_1(0))=(\pi,\ph_1^0)$ and 
$$\ue{F}(q(\delta),\ph_1(\delta),\ph_2^0,\delta)=\ue{0},\qquad \forall \delta
\in (-\bar{\delta}_1,\bar{\delta}_1) \mbox{,}$$ so 
$(q(\delta),\ph_1(\delta),\ph_2^0)$ is a transversal hetero/homoclinic point. Note that, as all the entries of $M$ are $O(\mu)$, then $\bar{\delta}$ is $O(\mu)$.\\
Defined $D:=-[\pl_q \Delta A_1 \pl_{\ph_1} \Delta A_2-\pl_q \Delta
A_2 \pl_{\ph_1} \Delta A_1]^{-1}$, we get 
$$
\de \frac{d(q,\ph_1)}{d \delta}(\delta) =
\mu^{-1}\left(%
\begin{array}{cc}
  \pl_q \Delta A_1 &  \pl_{\ph_1} \Delta A_1 \\
  \pl_q \Delta A_2 &  \pl_{\ph_1} \Delta A_2 \\
\end{array}%
\right)^{-1}
\left(%
\begin{array}{c}
  \pl_{\delta} F_1 \\
  \pl_{\delta} F_2 \\
\end{array}%
\right) =\mu^{-1}D
\left(%
\begin{array}{c}
  -\pl_{\ph_1} \Delta p \\
  \pl_{q} \Delta p \\
\end{array}%
\right)
$$
In which we have used (\ref{eq:deltapneqzero}b) and formulae (\ref{eq:isoenqphi}). Property
(\ref{eq:variabilitaphiuno}b) is a direct consequence.\\
By (\ref{eq:deltapneqzero}a) and by continuity, there exists 
$\bar{\delta}_2>0$ such that $\pl_q \Delta p \neq
0$ for all $\delta \in (-\bar{\delta}_2,\bar{\delta}_2)$.
By choosing $\bar{\delta}:=\min\{\bar{\delta}_1,\bar{\delta}_2\}$,
(\ref{eq:variabilitaphiuno}a) immediately follows. 
\endproof

\subsection{Transition chains}
\begin{defn}
A transition chain of length $N$ for (\ref{eq:hamiltonian}) is a (finite) sequence of tori $\ml{T}_k:=\ml{T}(\ue{A}^k)$ such that
\[ W^u(\ml{T}_k) \pitchfork W^s(\ml{T}_{k+1}),\qquad \forall k=1,\ldots,N-1   \mbox{.}\]
\end{defn}
By (\ref{eq:isoenergy}) the sequence $\{\ue{A}^k\}$ is necessarily of the form $\ue{A}(y_k)$, $y_k \in \R$, where
\[ \ue{A}(y):=\ue{A}^0+\ue{\omega}^{\perp}y \mbox{.}\]
Set $\delta^k:=y^{k+1}-y^k$, by lemma \ref{lemma:splitting} every sequence $\ml{T}_k$ such that $|\delta^k| \leq \bar{\delta}$ for all $k=1,\ldots,N-1$ is a transition chain. 
\begin{defn} We call \textbf{equally spaced chain} (ESC) a transition chain $\ml{T}_k$ such that 
\[ \delta^k=\{-\hat{\delta},0,\hat{\delta}\} \]
for some fixed $\hat{\delta} \leq \frac{3}{4} \bar{\delta}$.
\end{defn}
Now consider the sequence of intervals
\[ E_k :=(y^k-C_1\mu,y^k+C_1 \mu) \mbox{,}\]
with $C_1$ suitably chosen in a way that $\ml{T}(\ue{A}(\tilde{y}^k))$ is still a transition chain\footnote{It is sufficient to choose $C_1 \leq \bar{\delta}/(16 \mu)$.} for all $\tilde{y}^k \in E_k$ and for all $k$.
\begin{defn}
Let $\ml{T}_k^*$ an ESC and let $E_k^*$ be its associate sequence of intervals. An \textbf{equivalent elastic chain} (EEC) to the given ESC, is a sequence $\ml{T}_k$ such that $y_k \in E_k^*$ for all $k$.
\end{defn}
In this setting we have the following 
\begin{satz}\label{th:drift}
Let $N=O(\mu^{-1})$ and sufficiently small $\mu$. For all ESC of length $N$ there exists an EEC of the same length and a set of initial data $\ml{U}$ such that the solution of (\ref{eq:hamiltonian}) starting in $\ml{U}$ shadows\footnote{we do not need here a precise definition of the concept of shadowing as it will be intrinsically given by using windows.} the EEC with a time
\begin{equation}\label{eq:transitiontime}
T=O\left(\frac{1}{\mu} \ln \frac{1}{\mu}\right) \mbox{.}
\end{equation}
\end{satz}
As an immediate consequence we get
\begin{cor}[Arnold's Diffusion]
The Arnold's Diffusion takes place in the system (\ref{eq:hamiltonian}) on a set of initial data $\ml{U}_d \subset \ml{U}$ with a time given by (\ref{eq:transitiontime})  
\end{cor}
\proof
Choose the EEC given by $\delta^k:=\bar{\delta}$ (or $-\bar{\delta}$) for all $k$ and apply theorem \ref{th:drift}. As $\delta_k=O(\mu)$, an $O(\mu^{-1})$ number of transitions leads to a drift of $O(1)$. 
\endproof
\begin{rem}
As the set of all ESC of length $N$ has cardinality equal to $3^N$, the drift theorem \ref{th:drift} shows that the system (\ref{eq:hamiltonian}) possesses a rich behaviour. In particular, every step on the transition chain is conjugate to a random walk on a ternary tree. Arnold's diffusion is simply associated to the two lateral branches of the entire tree.
\end{rem}

\section{Flow approximation and windows}
This section is devoted to fix the background material in order to prove the theorem. First of all we are going to construct a suitable approximation (in the ``linear map + remainder'' form) of the Hamiltonian flow close to a generic sequence of hetero/homoclinic points associated with a given transition chain. This is obtained in a standard way, that is, by a composition of two maps: the first one describes the ``inner'' motion close to the invariant tori and the second one the ``outer'' motion from one torus to another. As a difference from \cite{mar96} and related works, we use Gallavotti's normal form instead of the Graff-Treshev one, giving a different criterion to determine the Jacobian matrix of the outer map. We also recall some basic notion on windows and their main properties for the reader convenience. For a more efficient comparison with our result, we shall follow notations of \cite{mar96} and of \cite{cresguil3} very closely from now on.

\subsection{Inner dynamics: normal form} 
Let $\ml{T}(\ue{A}^*)$ be an invariant torus for the system 
(\ref{eq:hamiltonian}). Fixed $\tilde{R},\tilde{\kappa}$, let us consider the phase space region defined as
$$
\tilde{U}_{\tilde{R},\tilde{\kappa}}(\ue{A}^*):=\{(q,\ue{\ph},p,\ue{A}):|p|,|q|<\tilde{R}, \quad
|\ue{A}-\ue{A}^*|<\mu \tilde{\kappa}, \quad \ue{\ph} \in
\mathbb{T}^2 \} \mbox{,}
$$ 
that can be regarded as a ``neighbourhood'' of $\ml{T}(\ue{A}^*)$.
\begin{satz}[Gallavotti, $1997$]\label{teo:formanormale}
For all $\ue{A}^*$ and $\mu$ small enough, there exist $\tilde{R},\tilde{\kappa},R,\kappa$ independent on $\mu$ and a canonical map
$$
\mathcal{C}: (q,\ue{\ph},p,\ue{A})\in
\tilde{U}_{\tilde{R},\tilde{\kappa}}(\ue{A}^*)  
\rw (Q,\ue{\psi},P,\ue{I}) \in U_{R,\kappa}(\ue{I}^*)
$$
with 
$$U_{R,\kappa}(\ue{I}^*):=\{(Q,\ue{\psi},P,\ue{I}):|P|,|Q|<R,\quad
|\ue{I}-\ue{I}^*|<\mu \kappa, \quad \ue{\psi} \in \mathbb{T}^2 \}
\mbox{,}$$ 
and $(0,\ue{\psi},0,\ue{I}^*)=\ml{C}(0,\ue{\ph},0,\ue{A}^*)$, casting the Hamiltonian (\ref{eq:hamiltonian}) in the following form
\footnote{this result holds with a more general perturbing function, in particular, it does not require the presence of the term $1-\cos q$.}
\begin{equation}\label{eq:hamnormali}
H(Q,\ue{\psi},P,\ue{I})=\ue{\omega} \cdot \ue{I}+J_{\mu}(\eta) \mbox{.}
\end{equation}
with $\eta:=PQ$ and where $J_{\mu}(\eta)$ is an analytic function on $U_{R,\kappa} \x (0,\mu_0)$ for some $\mu_0>0$.
\end{satz}
In the new variables set, the system is clearly integrable and its evolution is given by
\begin{equation}\label{eq:motouk}
 (Q(t),P(t),\ue{\psi}(t),\ue{I}(t))
=(Q(t_0)e^{g(\eta)t},P(t_0)e^{-g(\eta)t},\ue{\psi}(t_0)+\ue{\omega}(t-t_0),\ue{I}(t_0)) 
\end{equation}
where $g(\eta):=\frac{d}{d\eta}J_{\mu}(\eta)$ and $g(0)=:g$.\\
The proof of this result can be found in \cite{gal97}, and goes along the lines of the (more general) anisochronous case proof of \cite{chigal94}. In order to clarify some consideration we shall need later, let us recall that the canonical map $\ml{C}$ is the composition of two maps: the Jacobi map $\ml{C}_J$
\begin{equation}\label{eq:trasfj}
\left\{
\begin{array}{rcl}
q_J&=&S(q,p)\\
p_J&=&R(q,p)
\end{array}
\right.
\end{equation}
(straightening the invariant manifolds  of the unperturbed pendulum) and a KAM map $\ml{C}_K$ ($\mu$ close to the identity) 
\begin{equation}\label{eq:trasfkam}
\left\{
\begin{array}{rcl}
q_{J} & = & Q+\mu L(Q,\ue{\psi},P,\mu)\\
\ue{\ph} & = & \ue{\psi}\\
p_{J} & = & P+\mu M(Q,\ue{\psi},P,\mu)\\
\ue{A} & = & \ue{I}+\mu \ue{N}(Q,\ue{\psi},P,\mu)
\end{array}
\right.
\end{equation}
whose effect is to remove the perturbation at all orders.\\
\subsection{Reduced system of coordinates} 
In the neighbourhood $U_{R,\kappa}(\ue{I}^*)$, is useful to define a reduced system of coordinates. This standard approach is systematically used in such cases (see, e.g. in \cite{east81} or  \cite{mar96}) and can be obtained as follows 
\begin{enumerate}
\item Fix a Poincar\'{e} section on $\T^2$: due to isochrony, every straight line on $\T^2$ whose direction is not aligned with $\ue{\omega}$, is a transversal section. The simplest choice is 
$$\Sigma_N:=\{\ue{\psi} \in \T^2: \quad \psi_2=\ph_2^0\} \mbox{,}$$
for a complete interface with the ``external'' section $\Sigma$. In such a way we get a solutions discretization, 
and the sampling time $t^*$ is constant because of isochrony. Without loss of generality, we shall suppose $t^*$ to be equal to $1$.\\
As $\ue{\omega}$ is a Diophantine vector, the recurrence frequency on $\Sigma_N$, denoted as $\nu$, is a $(C,\tau)-$Diophantine number.
\item Set an energy level
$\ml{M}(h):=\{(P,Q,\ue{\psi},\ue{I}): H=h\}$ and suppress the action $I_2$ by restricting to $\ml{M}(h)$. 
\end{enumerate}
By defining $\theta:=\ph_1$ and $\rho:=I_1$, the coordinate system of the restriction to $V^h(\rho^*):=U_{R,\kappa}(\ue{I}^*) \cap
\Sigma_N \cap \ml{M}(h)$ is simply $(Q,\theta,P,\rho)$.
By taking into account of (\ref{eq:motouk}), the images of the section map  $f:V^h(\rho^*) \rw V^h(\rho^*)$, easily write as 
\begin{equation}\label{eq:mappapoincare}
\underline{f}^n(Q_0,\theta_0,P_0,\rho_0):=(Q_0 e^{n g(\eta)},\theta_0+n \nu,P_0
e^{-n g(\eta)},\rho_0) \mbox{,}
\end{equation}
where $\eta=Q_0P_0$. If we write 
$$e^{\pm n g(\eta)}=e^{\pm n g(0)}\left[1+\left(\pm n \eta \tilde{g}(x)+\frac{1}{2}(n \eta \tilde{g}(\eta) )^2\pm \ldots\right)\right]=:L^{\pm n}[1+F_{\pm}(\eta)] \mbox{,}$$ 
the map (\ref{eq:mappapoincare}) can be expanded in a neighbourhood of $\ue{X}_k$ in the following way
\begin{equation}\label{eq:expintmap}
 \ue{f}^n (\ue{X}_k+\underline{\xi})=(L^n(Q_k+\xi_1) ,\theta_k+\xi_2+n \nu,L^{-n} \xi_3, \rho_k+\xi_4 )+\underline{R}(\ue{X}_{k},\ue{\xi})\mbox{,}
\end{equation}
with 
\begin{equation}\label{eq:restof}
\ue{R}(\ue{X}_{k},\ue{\xi})=(L^n(Q_{k+1}+\xi_1) F_-(\eta),0, L^{-n} \xi_3 F_{+}(\eta),0) \mbox{,}
\end{equation}
where in such case $\eta=(Q_{k+1}+\xi_1)\xi_3$. We also denote with $\hat{\ue{f}}^n$ the linear part $\ue{f}^n-\underline{R}$.
\subsection{Sequences of hetero/homoclinic points associated to ESCs and EECs}\label{subsec:sequences}
By lemma \ref{lemma:splitting} we can associate with an ESC $\ml{T}_k$ a sequence of hetero/homoclinic points 
$(q^k,\ph_1^k):=(q(\delta_k),\ph_1(\delta_k))$ on $\Sigma$.\\
We denote with $\Phi^l$, $l \in \Z$, the sampling of the phase flow on  $\Sigma$.\\
Let $k$ be arbitrarily chosen and let
$$\ue{z}_k:=(q^k,\ph_1^k,\ph_2^0,p_0(q^k)+\mu
p^u(q^k,\ph_1^k,\ph_2^0,\mu),\ue{A}(y_k)+\mu
\ue{A}^u(q^k,\ph_1^k,\ph_2^0,\mu)) \mbox{.}$$ By construction 
$\ue{z}_k \in \ml{M} \cap \ml{M}(h) \cap \Sigma$.\\As $\lim_{l \rw -\infty} \Phi^l (\ue{z}_k) = \ml{T}_{k}$ and
$\lim_{l \rw +\infty} \Phi^l (\ue{z}_k) = \ml{T}_{k+1}$, for a sufficiently large $|l|$
we get $\Phi^{-l} (\ue{z}_k) \in
\tilde{U}_k$ and $\Phi^{l} (\ue{z}_k) \in \tilde{U}_{k+1}$. In this way, these points can be written in the respective normal coordinates of the tori. In particular, there exist $l_k^{\pm} \in \N$ such that\footnote{we denote e.g. with $(\ue{v})_{v_1}$ the $v_1$ component of $\ue{v}$.}
$$
(\ml{C} \circ \Phi^{-l_k^-}(\ue{z}_k))_Q<R/2,\qquad (\ml{C} \circ \Phi^{l_k^+} (\ue{z}_k))_P<R/2 \mbox{,}
$$
In the reduced system these points reads as
\begin{equation}\label{eq:puntieteroclìni}
\begin{array}{rcl}
\ue{X}_k=(Q_k,\theta_k,0,\rho_k )&:=&(\ml{C}\circ \Phi^{l_k^-}(\ue{z}_k))|_{V_k}\\
\ue{X}_k'=(0,\theta_k',P_k',\rho_{k+1})&:=&(\ml{C}\circ
\Phi^{l_k^+}(\ue{z}_k))|_{V_{k+1}}
\end{array}\mbox{.}
\end{equation}
Let us consider an ESC, for all EECs associated with the given ESC, by following the previous construction we get a double sequence of hetero/homoclinic points $\{X_k(\tilde{y}_k),X_k'(\tilde{y}_k)\}_{k=1,\ldots,N}$ for all $\tilde{y}_k \in E_k$. If necessary, we could reduce the range of variation $\bar{\delta}$ (i.e. $C_1$) in such a way $Q_k(y),P_k'(y) <(3/4)R$ for all $y \in E_k$.\\The key point is that, due to isochrony, by varying $y$ in $E_k$ (of an $O(\mu)$), we can move the coordinates of the hetero/homoclinic point (of an $O(1)$, see (\ref{eq:variabilitaphiuno}b)) on $\T^2$ without ever falling out of $\Sigma_N$, i.e. by functions 
\begin{equation}\label{eq:thetaky}
\theta_k(y):=\ph_1(\delta)-l_k^- \nu, \qquad \theta_k'(y):=\ph_1(\delta)+l_k^+ \nu \mbox{,}
\end{equation}
where $\delta:=y_{k+1}-y$ with $y_{k+1}$ fixed in $E_{k+1}$.\\
This is the advantage in our formulation of lemma \ref{lemma:splitting} with respect to the setting used in \cite{gal99}. As $\Sigma_N$ is one dimensional the problem of finding $y$ such that e.g. $\theta_k(y)$ is equal to some given $\theta^*$ in the image of $\theta_k(y)$ is achievable in a constructive way\footnote{compare with comments in \cite[pag. $307$ and $309$]{gal99}.}.\\ 
As a direct consequence of the normal form, the manifolds $W^u(\ml{T}_{k})$ and $W^s(\ml{T}_{k+1})$ in $V_k$ and $V_{k+1}$ respectively, are given by 
\[ W^u(\ml{T}_{k})=\{(Q,\theta,P,\rho):P=0,\,\rho=\rho_k\},
\qquad
W^{s}(\ml{T}_{k+1})=\{(Q,\theta,P,\rho):Q=0,\,\rho=\rho_{k+1}\}
 \mbox{.}\]
As $l_k^-=O(1)$ (see remark \ref{rem:splitting}), by writing down the expression of the whisker $W^s(\ml{T}_{k+1})$ in terms of normal coordinates\footnote{so backward evolved by a time $l_k^-$.}, it is possible to recover (up to higher orders) two smooth functions $P_{k+1}(Q,\theta;y)$ and $\rho_{k+1}(Q,\theta;y)$ such that the connected component of $W^s(\ml{T}_{k+1})$ intersecting $W^u(\ml{T}_k)$ at $\ue{X}_k(y)$ in $V_k$ is parameterized by
\begin{equation}\label{eq:ykpuno}
\ue{Y}_{k+1}(Q,\theta;y):=\{(Q,\theta,P_{k+1}(Q,\theta;y),\rho_{k+1}(Q,\theta;y))\}, 
\qquad 
\ue{Y}_{k+1}(y)(Q_k,\theta_k)=\ue{X}_k(y) \mbox{.}	
\end{equation}
Furthermore, the normal splitting matrix
\begin{equation}\label{eq:matricesplittingnormali}
    M_N:=\left(%
\begin{array}{cc}
  \pl_Q P_{k+1} & \pl_{\theta} P_{k+1} \\
  \pl_Q \rho_{k+1} & \pl_{\theta} \rho_{k+1} \\
\end{array}%
\right)(Q_k,\theta_k;y)
\end{equation}
satisfies
\begin{equation}\label{eq:quellochesoddisfasplittingnormali}
(M_N)_{11} \neq 0,\qquad rk M_N=2\mbox{,}
\end{equation}
by virtue of (\ref{eq:deltapneqzero}) for all $y \in E_k$ (if necessary we could reduce $C_1$ further). As a geometrical consequence, we have that $W^s(\ml{T}_{k+1})$ is still a graph over the angles in the normal (reduced) coordinates system.

\subsection{Outer dynamics: transition map}
As we have written by $\ue{f}^n$ the evolution in $V_k$, we want to give, here, an approximation of the Hamiltonian flow from $V_k$ to $V_{k+1}$. For this purpose, let $\underline{\Psi}:V_k \rw V_{k+1}$ the Hamiltonian evolution implicitly defined by 
\begin{equation}\label{psi}
\underline{\Psi}(\ue{X}_k)=\ue{X}_k' \mbox{.}
\end{equation}
By regular dependence on initial data, there exist neighbourhoods 
$\ml{D}_k$ of $\underline{X}_k$ and $\ml{D}_k'$ of $\underline{X}_k'$ respectively, such that 
$\ue{\Psi}$ maps $\ml{D}_k$ into $\ml{D}_k'$ and, for all $\ue{X}_k+\ue{\xi} \in \ml{D}_k$ can be written in the form   
\begin{equation}\label{eq:mappapsi}
\ue{\Psi}(\ue{X}_k+\ue{\xi})=\ue{X}_k'+ A \ue{\xi}  + \ue{\tilde{R}}(\ue{X}_k,\ue{\xi}) \mbox{.}
\end{equation}
We have denoted by $A:=\mathcal{D}\ue{\Psi}(\ue{X}_k)$ the Jacobian matrix of the map computed at the hetero/homoclinic point and by $\ue{\tilde{R}}(\ue{X}_k,\ue{\xi})$ the remainder that in terms of components reads as
\begin{equation}\label{eq:restopsi}
\tilde{R}^i(\ue{X}_k,\ue{\xi})=\int_0^1 (1-t) \ue{\xi} \ml{H} \Psi^i(\ue{X}_k+t
\ue{\xi}) \cdot \ue{\xi}^T dt,\qquad i=1,\ldots,4
\end{equation}
$\ml{H}$ denotes the Hessian matrix operator. We denote with $\hat{\ue{\Psi}}$ the linear part $\ue{\Psi}-\ue{\tilde{R}}$.\\
The map $\ue{\Psi}$, also known as \emph{separatrix map}, see e.g. \cite{tre}, is a standard and widely used tool to construct trajectories close to invariant manifolds. The approach we are going to use, in order to describe the structure of $A$, is of geometrical nature and it seems to be very close to the one used in \cite{tre}. The matrix $A$ is essential for a non-trivial approximation of the flow, however it is not  difficult to understand that this matrix conceals a complicate structure in our case. For instance,  let us denote with $\Phi_{k}'$ the Hamiltonian flow (in global coordinates) mapping $\ml{D}_k$ in $\ml{D}_k'$ and $\ml{C}_R^{k,k+1}$ the canonical map $\ml{C}$ with the coordinates reduction to $V_k$ and $V_{k+1}$ respectively. So, by construction we have  
$$\ue{\Psi}(Q,\theta,P,\rho)=(\mathcal{C}_R^{k+1} \circ \Phi_{k}' \circ (\mathcal{C}_R^k)^{-1})(Q,\theta,P,\rho) \mbox{.}$$
As one can easily deduce by (\ref{eq:trasfj}) and (\ref{eq:trasfkam}), this computation turns out to be a very difficult task\footnote{the main difficulty arises if one write down the perturbing function in terms of Jacobi's coordinates and then attempt to perform a perturbative step.}, even at the first order in $\mu$, unless one would like to resort to a numerical approach.\\ 
We state the following 
\begin{prop}
The matrix $A$ for the system (\ref{eq:hamiltonian}) takes the form
\begin{equation}\label{eq:formageneralea}
A=
\left(
\begin{array}{cccc}
0 & 0 & a & 0\\
0 & 1 & 0 & 0\\
-1/a & 0 & 0 & 0\\
0 & 0 & 0 & 1
\end{array}
\right)
+
\mu
\left(
\begin{array}{cccc}
a_{11} & a_{12} & a_{13} & 0\\
0 & 0 & 0 & 0\\
a_{31} & a_{32} & a_{33} & 0\\
a_{41} & a_{42} & a_{43} & 0
\end{array}
\right) \mbox{,}
\end{equation}
with $a \neq 0$, and\footnote{In spite of the dependence also on $k$ (as $A$ is the Jacobian matrix computed at $\ue{X}_k$), we suppress the index $k$ to avoid a cumbersome notation.} $a_{ij}:=a_{ij}^{k}(\mu)$. These satisfies, for sufficiently small $\mu$ 
\begin{equation}\label{eq:splittingaij}
a_{11} \neq 0, \qquad D:=a_{11}a_{42}-a_{41}a_{12} \neq 0 \mbox{.}
\end{equation}
\end{prop}
\noindent Note that the above structure holds at all orders in $\mu$. We denote with $A_0$ the unperturbed part of $A$.
\proof
If $\mu=0$ the system is integrable, then rotators and pendulum are uncoupled. Moreover, the canonical map $\mathcal{C}$ reduces to the Jacobi's map for the pendulum, whose Hamiltonian takes simply the form $J_0(PQ)$.\\
Let $\Pi$ the projection operator on the $(P,Q)$ plane. By conservation of energy, for $(P,Q) \in \Pi V_k$ and $(P',Q') \in \Pi V_{k+1}$ we get $J_0(P Q)=J_0(P' Q')$. Now  consider a small increment of $(Q_k,0)$ of the form $(Q_k+\sigma,\delta)$ by denoting with $(f(\sigma,\delta),P_k'+g(\sigma,\delta))$ corresponding variations in $\Pi V_{k+1}$. It is understood that $f(\delta,\sigma)$ and $g(\delta,\sigma)$ are regular functions by definition of $\ue{\Psi}$ and such that $f(0,0)=g(0,0)=0$. Denote with $u_{ij}$ the entries of $\frac{\pl(f,g)}{\pl(\sigma,\delta)}(0,0)$. By conservation of energy, turns out to be, up to higher orders in $(\delta,\sigma)$, 
$$u_{11} u_{21}=0,\quad u_{12} u_{22}=0, \quad u_{11}u_{22}+u_{21}u_{12}=1 \mbox{.}$$
Now recall that the intersection between $W_k^u$ and  $W_{k+1}^s$ is flat, then increments in direction $(Q,0)$ in $\Pi V_k$ imply displacements in direction $(0,-P)$ in $\Pi V_{k+1}$, from which $u_{21}<0$. So, by the previous equalities we get $u_{11}=0$, so $u_{12}=1/u_{21}$ and then $u_{22}=0$. It is sufficient to define $a:=-u_{21}$ and note that the Jacobian matrix in the $(\theta,\rho)$ variables is trivially the identity matrix, to get the required form of the unperturbed $A$.\\
Now suppose $\mu \neq 0$. By definition, if $\ue{\Psi}$ is restricted on the graphs of the functions $P_{k+1}(Q,\theta),\rho_{k+1}(Q,\theta)$
(parameterizing $W_{k+1}^s$ in $V_k$), its image has to produce points on $W_{k+1}^s$ in $V_{k+1}$ and then such that $Q=0$
and $\rho=\rho_{k+1}$. In other terms, for sufficiently small there exists $\bar{r}>0$ such that,
for all $(\sigma,\ep) \in \ml{B}_{\bar{r}}(\underline{0})$

\begin{equation}\label{eq:annullamentopsi}
\left\{
\begin{array}{rcl}
    \Psi_1(Q_k+\sigma, \theta_k+\ep, P_{k+1}(Q_k+\sigma,\theta_k+\ep),\rho_{k+1}(Q_k+\sigma,\theta_k+\ep)) & = & 0\\
    \Psi_4(Q_k+\sigma, \theta_k+\ep, P_{k+1}(Q_k+\sigma,\theta_k+\ep),\rho_{k+1}(Q_k+\sigma,\theta_k+\ep)) & = & \rho_{k+1}
\end{array}
\right. \mbox{.}
\end{equation}
\nt By expanding the first equation in a neighbourhood of $(\sigma,\ep)=\underline{0}$, we get
\begin{equation}\label{eq:duerelazioni}
\left\{
\begin{array}{rcl}
a_{11}+a_{13}\frac{\pl P_{k+1}}{\pl Q}+a_{14}\frac{\pl \rho_{k+1}}{\pl Q}&=&0\\
a_{12}+a_{13}\frac{\pl P_{k+1}}{\pl \theta}+a_{14}\frac{\pl \rho_{k+1}}{\pl \theta}&=&0
\end{array}
\right. \mbox{,}
\end{equation}
where we have recognized the elements of $A$. By the unperturbed form of $A$ we know $a_{1j}=\mu
a_{1j}^1+o(\mu)$ for $j=1,2,4$ while $a_{13}=a+\mu
a_{13}^1+o(\mu)$. On the other hand, the derivatives appearing in the previous equation are exactly those of the splitting matrix $M_N$ (\ref{eq:matricesplittingnormali}). Then, by denoting with 
$$\mu
\left(
\begin{array}{cc}
h_{11} & h_{12}\\
h_{21} & h_{22}
\end{array}
\right):=
\left(
\begin{array}{cc}
\frac{\pl P_{k+1}}{\pl Q} & \frac{\pl P_{k+1}}{\pl \theta}\\
\frac{\pl \rho_{k+1}}{\pl Q} & \frac{\pl \rho_{k+1}}{\pl \theta}
\end{array}
\right) \mbox{,}
$$
and equating powers of $\mu$ in (\ref{eq:duerelazioni}) we get $a_{11}^1=-a h_{11}$ and $a_{12}^1=-a h_{21}$. First of all we note that $a_{11} \neq 0$ for $\mu \neq 0$ by (\ref{eq:quellochesoddisfasplittingnormali}a). In similar way, the second equation of (\ref{eq:annullamentopsi}) yields $a_{41}^1=- h_{21}$ and $a_{42}^1=- h_{22}$. Then, by keeping in mind the definition of $D$, follows 
$$D=a(h_{11}h_{22}-h_{12}h_{21}) \mbox{,}$$
that is, up to higher orders, the determinant of $M_N$: the latter is non-zero by 
(\ref{eq:quellochesoddisfasplittingnormali}b).\\
Note that the argument is until now quite general, and can be reasonably extended to more general systems.\\
Now we attempt to get some additional information by using the very simple form of our system, in order to simplify, as much as possible, the structure of the perturbed part of $A$. \\
Let us consider two matrices $B_1,B_2 \in \mathbb{M}(n,m)$ for some $n,m \in \N$. We say that $B_1$ and $B_2$ possess the \emph{same structure}, and we denote this property with $B_1 \approx B_2$, if $B_1$ and $B_2$ have null entries exactly in the same position(s).\\
So, by (\ref{eq:trasfj}) and (\ref{eq:trasfkam}) the Jacobian matrices of the canonical maps satisfy 
$$
D \mathcal{C}_J  \equiv 
\left(
\begin{array}{cccccc}
S_q & 0 & 0 & S_p & 0 & 0\\
0 & 1 & 0 & 0 & 0 & 0\\
0 & 0 & 1 & 0 & 0 & 0\\
R_q & 0 & 0 & R_p & 0 & 0\\
0 & 0 & 0 & 0 & 1 & 0\\
0 & 0 & 0 & 0 & 0 & 1
\end{array}
\right),\qquad 
D \mathcal{C}_K  \equiv  \mathbb{I}+\mu
\left(
\begin{array}{cccccc}
L_Q & L_{\psi_1} & L_{\psi_2} & L_P & 0 & 0\\
0 & 0 & 0 & 0 & 0 & 0\\
0 & 0 & 0 & 0 & 0 & 0\\
M_Q & M_{\psi_1} & M_{\psi_2} & M_P & 0 & 0\\
N_Q^1 & N_{\psi_1}^1 & N_{\psi_2}^1 & N_P^1 & 0 & 0\\
N_Q^2 & N_{\psi_1}^2 & N_{\psi_2}^2 & N_P^2 & 0 & 0\\
\end{array}
\right) \mbox{,}
$$
where $\mathbb{I}$ is the identity matrix. Both matrices are computed, by definition, at the hetero/homoclinic point and the entries of the second one are functions of $\mu$.\\
First of all is clear that $D\mathcal{C}_J \approx D \mathcal{C}_J^{-1}$. Furthermore, it is easy to check from (\ref{eq:trasfkam}) that $D\mathcal{C}_K \approx D \mathcal{C}_K^{-1}$.\\
Now, simply by looking at the equations of motion, we can see that the pendulum evolution does not depends on the action variables at all orders in $\mu$. This implies that the actions evolution (also due to the independence of the perturbing function on $\underline{A}$) is independent on $\underline{A}$ at $O(\mu^k)$ for all $k \geq 1$. So the Jacobian matrix of the flow satisfies 
$$
D \Phi^t \approx
\left(
\begin{array}{cccccc}
\x & 0 & 0 & \x & 0 & 0\\
0 & 1 & 0 & 0 & 0 & 0\\
0 & 0 & 1 & 0 & 0 & 0\\
\x & 0 & 0 & \x & 0 & 0\\
0 & 0 & 0 & 0 & 1 & 0\\
0 & 0 & 0 & 0 & 0 & 1
\end{array}
\right)
+\mu
\left(
\begin{array}{cccccc}
\x & \x & \x & \x & 0 & 0\\
0 & 0 & 0 & 0 & 0 & 0\\
0 & 0 & 0 & 0 & 0 & 0\\
\x & \x & \x & \x & 0 & 0\\
\x & \x & \x & \x & 0 & 0\\
\x & \x & \x & \x & 0 & 0
\end{array}
\right) \mbox{,}
$$
where the entries of the second matrix are functions of $\mu$.\\
In conclusion, by definition, we have that $A := D \ml{C}_K \cdot D \ml{C}_J \cdot D \Phi^t \cdot D \ml{C}_J^{-1} \cdot D \ml{C}_K^{-1} \approx D \Phi^t$. By reducing to dimension $4$ we get the perturbed part of (\ref{eq:formageneralea}) .
\endproof
\begin{rem}
Clearly, the previous argument does not guarantee that the entries of $A$ marked by $\x$ are non-zero since some of them could vanish due to some symmetry of the product. As we shall see, the only key properties we shall need are given by  (\ref{eq:splittingaij}): it does not matter if the other entries vanish or not. In this sense we agree with the topological criterion used in \cite[p. 239]{mar96} to give the structure of his matrix $A$.
\end{rem}

\subsection{Windows}
\begin{defn}
Let $M \subset \R^d$, and $d_h,d_v \in \N$ such that $d_h+d_v=d$. let $L_h:=[-1,1]^{d_h}$, $L_v:=[-1,1]^{d_v}$ and $\ml{Q}:=L_h \times L_v$.\\
A \tf{window} is a $C^1$ diffeomorphism 
$$\ue{W}:\ml{Q} \rw M \mbox{.}$$
We denote with $\tilde{W}$ the image\footnote{we shall refer as widow either the map or its image.} of $\ml{Q}$ through $\ue{W}$.\\
The \emph{horizontals} of $\ue{W}$ are defined as the images of $\ue{W}(x_h,x_v)$ where $x_v \in L_v$ is kept fixed and $x_h$ varies in $L_h$. Similarly, the \emph{verticals} of $\ue{W}$ are defined as the images of $\ue{W}(x_h,x_v)$ with $x_h \in L_h$ fixed and variable $x_v \in L_v$.
\end{defn}
The point $\ue{W}(\underline{0})$ is called \emph{centre} of the window. We shall refer to \emph{non-degenerate} windows as windows containing a non-zero volume. 
Even if a non-necessary condition, it is natural to think the integers $d_h$ and $d_v$ as satisfying $d_h=d_v=d/2$. Actually, in the applications, $d$ is the (even) dimension of the (eventually reduced) phase space.\\ 
\nt An interesting case is represented by \emph{affine windows}. This class of windows arises when the function $\ue{W}$ is linear, and then can be written in the form 
$$\ue{W}^a(\ue{x})=\ue{c}+W \ue{x} \mbox{.}$$
The point $\ue{c}$ is the centre and $W$ is a $d \x d$ matrix we shall call \eh{representative matrix} of $\ue{W}^a$. Clearly, a window is non-degenerate if $\det W \neq 0$. 
\begin{defn}\label{def:correttoall}
We say that a window $\ue{W}_1$ is \tf{correctly aligned on} (c.a.o.)
$\ue{W}_2$ if for all $y_v \in L_v$ and for all $y_h \in L_h$ there exists unique $x_h \in (-1,1)^{d_h}$ and $x_v \in (-1,1)^{d_v}$ such that 
\begin{equation}\label{eq:correttoallineamento}
\ue{W}_1(x_h,y_v)=\ue{W}_2(y_h,x_v) \mbox{,}
\end{equation}
and the intersection is transversal.
\end{defn}
\nt The notion of correct alignment leads to the following key result
\begin{satz}[Shadowing, Easton]
Let $\{\ml{D}_k\}_{k=1,\ldots,N}$ be a collection of open sets of $\R^d$, $\Psi_k : \ml{D}_{k} \rw \ml{D}_{k+1}$ a family of diffeomorphisms and  $\ue{W}_k$ a family of windows $\ue{W}_k:\ml{Q} \rw \ml{D}_k$. Let, for all $k$, $\ue{W}_{k}':=\Psi_k \circ  \ue{W}_k$.\\
If $\ue{W}_k'$ c.a.o.  $\ue{W}_{k+1}$ for all $k=1,\ldots,N-1$ there exists (at least) a point $\ue{x}^*$ of $\tilde{W}_1 \subset \ml{D}_1$ such that $\Psi_{N-1} \circ \ldots \circ \Psi_1 (\ue{x}^*) \in \tilde{W}_N \subset \ml{D}_N$.
\end{satz}
\noindent For the proof, given in a more general context, we refer to \cite{eastmcg79}.\\

\subsection{A correct alignment criterion}
The following statement, due to J.P. Marco, give us a set of sufficient conditions in order to establish if a window is correctly aligned on another one. First of all, some notations:\\
Let $\mathbb{M}(n,m):=GL(\R,n \times m)$ be the linear group of the $n \x m$ matrices with real entries. For all $\ue{v} \in \R^d$ let $\norm{\ue{v}}_{\infty}:=\max_{j=1,\ldots,d}|v_j|$. Consequently, for all $A \in \mathbb{M}(n,m)$ we denote $\norm{A}_{\infty}:=\max_{i=1,\ldots,n} \sum_{j=1}^m |a_{ij}|$.\\
Now, let $\ue{F}(\ue{x}) \in C^k(\ml{U},\R^n)$, $k \geq 1$. The standard $C^1$ norm over $\ml{U}$ is denoted as 
\begin{equation}
\norm{\ue{F}(\ue{x})}_{C^1(\ml{U})}:=\max\{\sup_{\ue{x} \in \ml{U}}\norm{\ue{F}(\ue{x})}_{\infty},\sup_{\ue{x} \in \ml{U}} \norm{D \ue{F}(\ue{x})}_{\infty}\} \mbox{.}	
\end{equation}
\begin{lem}[Marco, $1996$]\label{lem:marco}
Let
\[
\underline{W}_{1}^{a}=\underline{c}_{1}+W_{1}\underline{x},\qquad\underline{W}_{2}^{a}=\underline{c}_{2}+W_{2}\underline{x} \mbox{,}
\]
be two affine windows, where $\underline{c}_{1,2}\in\mathbb{R}^{d}$,
and $W_{1,2}\in\mathbb{M}(d,d)$ are of the form 
\[
W_{i}=\left(\begin{array}{cc}
W_{i}^{1} & W_{i}^{3}\\
W_{i}^{2} & W_{i}^{4}
\end{array}\right) \mbox{,}
\]
with $W_{i}^{1}\in\mathbb{M}(d_{h},d_{h})$, $W_{i}^{2}\in\mathbb{M}(d_{h},d_{v})$,
$W_{i}^{3}\in\mathbb{M}(d_{v},d_{h})$ and $W_{i}^{4}\in\mathbb{M}(d_{v},d_{v})$.\\
Now define the {}``intermediary matrices''
\[
M=\left(\begin{array}{cc}
W_{1}^{1} & -W_{2}^{3}\\
W_{1}^{2} & -W_{2}^{4}
\end{array}\right),\qquad N=\left(\begin{array}{cc}
-W_{2}^{1} & W_{1}^{3}\\
-W_{2}^{2} & W_{1}^{4}
\end{array}\right) \mbox{,}
\]
and 
\begin{equation}\label{eq:condallineamentoaffine}
\chi^a:=\sup_{\underline{x} \in \mathcal{Q}} \norm{M^{-1}N \underline{x}+M^{-1}(\underline{c}_{2}-\underline{c}_{1})}_{\infty}\mbox{.}
\end{equation}
Suppose that $M$ is invertible and $\chi^a<1$. Then $\underline{W}_{1}^{a}$ c.a.o. $\underline{W}_{2}^{a}$.\\
Now consider the following windows 
\begin{equation}\label{eq:problemacompleto}
\underline{W}_{1}(\underline{x})=\underline{W}_{1}^{a}(\underline{x})+\underline{\hat{W}}_{1}(\underline{x}),\qquad
\underline{W}_{2}(\underline{x})=\underline{W}_{2}^{a}(\underline{x})+\underline{\hat{W}}_{2}(\underline{x}) \mbox{,}
\end{equation}
where the maps $\underline{\hat{W}}_{i}(\underline{x}) \in C^1(\ml{U}
,\R^d)$, and $\ml{U}$ is an open set containing $ \ml{Q}$. And define
$$\chi^c:=||M^{-1}(\underline{\hat{W}}_2(\underline{x})-\underline{\hat{W}}_1(\underline{x}))||_{C^1(\ml{U})}\mbox{,}$$  
Suppose that
\begin{enumerate}
\item[c1.] $\chi^a<1$, implying that $\underline{W}_{1}^{a}$ c.a.o. $\underline{W}_{2}^{a}$, 
\item[c2.]
\begin{equation}\label{eq:condizionechic}
    \chi^c<\frac{1}{4},\qquad \chi^a+\frac{\chi^c}{1-\chi^c}<1 \mbox{.}
\end{equation}
\end{enumerate}
Then $\underline{W}_{1}$ c.a.o.
$\underline{W}_{2}$.
\end{lem}
\nt As in \cite{mar96} it is natural to choose $\ml{U}:=(-2,2)^d$. A detailed proof of lemma \ref{lem:marco} can be found in the already mentioned paper or, more closely to the above formulation, in \cite{locmar05}. 

\subsection{A suitable class of windows}
Our aim is now to construct a prototype of non-degenerate affine windows for the problem at hand. Following \cite{east81} we choose a class of windows satisfying  
\begin{cond}\label{cond:geomcomp}
\begin{enumerate}
    \item The horizontals of $\underline{B}(\underline{x})$ are parallel to $T_{\underline{X}_k}W_k^u$,
    \item The verticals of $A \underline{B}(\underline{x})$ are parallel to\footnote{
    This item is an easier to handle condition than ``the verticals of $\ue{\Psi} \circ \ue{B}(\underline{x})$ are parallel to $T_{\underline{X}_k'}W_{k+1}^s$'', but the geometrical meaning is the same. In this way, since $T_{\underline{X}_k}W_k^u$ and $T_{\underline{X}_k}W_{k+1}^s$ are transversal by virtue of the splitting, horizontals and verticals of the window are transversal too, implying the non-degeneration property.} $T_{\underline{X}_k'}W_{k+1}^s$.
\end{enumerate}
\end{cond}
\noindent Before restricting ourselves to the particular form of $A$ given by (\ref{eq:formageneralea}), we give the following 
\begin{prop} Suppose $A=A_0+\mu A_1$ where $\{A_1\}_{ij}$ are generic functions of $\mu$ satisfying (\ref{eq:splittingaij}). There exist (at least) $\infty^4$ non-degenerate windows satisfying condition \ref{cond:geomcomp}. More precisely, let $b_{1,..,4} \in \mathbb{R} \setminus \{0\}$ be such parameters, the simplest\footnote{i.e. with the lowest number of non-vanishing entries, up to an entries interchanging.}  representative matrix takes the form 
\begin{equation}\label{eq:formasempliceb}
B=
\left(
\begin{array}{cccc}
b_1 & 0 & b_3 & 0\\
0 & b_2 & 0 & b_4\\
0 & 0 & \mu \Sigma_{33} b_3 & \mu \Sigma_{34} b_4\\
0 & 0 &  \mu \Sigma_{43} b_3 &  \mu \Sigma_{44} b_4
\end{array}
\right) \mbox{,}
\end{equation}
such that $\det B \neq 0$. $\Sigma_{kl}$ are suitable functions of $a_{ij}$.
\end{prop}
This structure allows us to simplify as much as possible the computation of the intermediary matrices without any loss of generality.\\
It will be useful to isolate a scaling factor in the parameters $b_i$. For this purpose we redefine $b_{i}=: \mu^{p} b_{i}$, and from now on $b_{i}$ are $O(1)$ constants and $p \in \mathbb{N}$ is a control parameter of the window size.\\ 
As in this circumstance, the particular structure of the perturbed part of $A$ given in (\ref{eq:formageneralea}) does not lead to substantial simplifications, so we prefer to bring for the moment a generic structure. In any case the condition $D \neq 0$ points out a clear relation between splitting and condition \ref{cond:geomcomp}.
\proof
Consider a generic $B=\{b_{ij}\}$, the above mentioned tangent spaces writes
$$
T_{\underline{X}_k}W_k^u=\{P=0,\quad \rho=\rho_k \},\qquad T_{\underline{X}_k'}W_{k+1}^s=\{Q=0,\quad \rho=\rho_{k+1}\} \mbox{.}
$$
Let us denote with $\ue{u}_i$ the canonical basis vectors of $\ml{Q}$ and with $\ue{e}_i$ those of $\R^4$. Then conditions \ref{cond:geomcomp} and non-degeneracy take the form
\begin{enumerate}
\item $    B \ue{u}_i \cdot \ue{e}_j  = 0   \qquad i=1,2 \quad j=3,4 $,
\item $   A B \ue{u}_i \cdot \ue{e}_j  = 0 \qquad i=3,4 \quad j=1,4 $,
\item $B$ invertible.
\end{enumerate}
The first condition is satisfied if $ b_{31}=b_{41}=b_{32}=b_{42}=0 $,
while the second one leads to the following linear systems 
$$
\left\{
\begin{array}{rcl}
f_1(\mu)b_{33}+ \mu a_{14} b_{43} & = & \mu d_1 \\
\mu a_{43} b_{33} +f_2(\mu) b_{43} & = & \mu d_3
\end{array}
\right.
\quad\mbox{;}\quad 
\left\{
\begin{array}{rcl}
f_1(\mu)b_{34}+ \mu a_{14} b_{44} & = & \mu d_2\\
\mu a_{43} b_{34} +f_2(\mu) b_{44} & = & \mu d_4
\end{array}
\right. \mbox{,}
$$
where 
$$
\begin{array}{rcl}
f_1 (\mu) &:= &(a+\mu a_{13})\\
f_2(\mu) &:= &1+ \mu a_{44}\\
\end{array}
\quad\mbox{;}\quad 
\begin{array}{rcl}
d_1 & := & -(a_{11} b_{13} + a_{12} b_{23})\\
d_3 & := & -(a_{41} b_{13} +  a_{42} b_{23})\\
d_2 & := & -(a_{11} b_{14} + a_{12} b_{24})\\
d_4 & := & -(a_{41} b_{14} +  a_{42} b_{24})
\end{array}
\mbox{.}
$$
By setting $\delta(\mu)=f_1(\mu) f_2(\mu)-\mu^2 a_{14}a_{43}=a+O(\mu)$ that is non-zero for sufficiently small $\mu$, the above systems give
\begin{equation}\label{eq:soluzb}
\begin{array}{rcl}
b_{33}&=& \mu \frac{d_1 f_2(\mu)-\mu a_{14}d_3}{\delta(\mu)}\\
b_{43}&=& \mu \frac{d_3 f_1(\mu)-\mu a_{43}d_1}{\delta(\mu)}
\end{array}
\quad\mbox{;}\quad 
\begin{array}{rcl}
b_{34}&=& \mu \frac{d_2 f_2(\mu)-\mu a_{14}d_4}{\delta(\mu)}\\
b_{44}&=& \mu \frac{d_4 f_1(\mu)-\mu a_{43}d_2}{\delta(\mu)}
\end{array}
\mbox{.}
\end{equation}
Note that, by the previous solution and by definition of $\delta(\mu)$, we get
\begin{equation}\label{eq:detfinestre}
b_{33}b_{44}-b_{43}b_{34}=\frac{\mu^2}{\delta(\mu)}D(b_{13}b_{24}-b_{23}b_{14}) \mbox{.}
\end{equation}
The simplest form of $B$ is then achievable by choosing
\begin{equation}\label{eq:semplb}
b_{23}=b_{14}=b_{21}=b_{12}=0  \mbox{,}
\end{equation}
in such a way 
\begin{equation}\label{eq:detb}
\det B=\frac{\mu^2}{\delta(\mu)} b_{11}b_{22}b_{13}b_{24} D\mbox{.}
\end{equation}
Now define $b_1:=b_{11}$, $b_2:=b_{22}$, $b_3:=b_{13}$ and $b_4:=b_{24}$. By hypothesis $D \neq 0$, so $B$ is invertible.\\
In conclusion, by substituting the simplified $d_{1,..,4}$ (by choices (\ref{eq:semplb}) in (\ref{eq:soluzb})) and defining
\begin{equation}\label{eq:sigmaij}
\begin{array}{rcl}
\Sigma_{33}&:=&\frac{-a_{11}+\mu (a_{14}a_{41}-a_{11}a_{44})}{\delta(\mu)}\\
\Sigma_{34}&:=&\frac{-a_{12}+\mu (a_{14}a_{42}-a_{12}a_{44})}{\delta(\mu)}
\end{array}
\quad\mbox{;}\quad 
\begin{array}{rcl}
\Sigma_{43}&:=&\frac{-a a_{41}+\mu (a_{11}a_{43}-a_{13}a_{41})}{\delta(\mu)}\\
\Sigma_{44}&:=&\frac{-a a_{42}+\mu (a_{12}a_{43}-a_{13}a_{42})}{\delta(\mu)}
\end{array}
\mbox{,}
\end{equation}
we get the required form. 
\endproof

\section{Proof of the theorem}
By using the preparatory material of the previous section, we want to prove theorem \ref{th:drift} via WM. The construction is standard: we consider the evolution through $\ue{\Psi}$ of a generic window $\ue{B}_k$ (``close'' to the $k-$th torus) obtaining a window $\ue{W}_1^k$ and the backward evolution through $\ue{f}^{-n}$ of a window $\ue{B}_{k+1}$ (close to the $k+1-$th torus), getting $\ue{W}_2^k$. Now we wonder if centres and representative matrices of $\ue{B}_{k,k+1}$ can be determined in such way $\ue{W}_1^k$ c.a.o. $\ue{W}_2^k$. The windows $\ue{B}_k$ possess, by construction, a ``shape'' that makes this property easier to obtain. Nevertheless, the absence of an anisochrony term, implies that the windows, if constructed as in \cite{mar96}, are not ``deformed'' in a suitable way to get correct alignment. This will be clear imposing the condition required by lemma \ref{lem:marco}, we use here as a guideline for the proof.
\subsection{Affine problem}
Let us consider a family of affine windows of the form 
$$\ue{B}_k(\ue{x})=\ue{p}_k+B_k \ue{x} \mbox{.}$$
defined in the neighbourhood of $\ue{X}_k$ 
\begin{equation}\label{eq:centriquattrod}
\ue{p}_k=(\sigma_k,0,\delta_k,0)+\ue{X}_k \mbox{,}
\end{equation}
with $B_k,\sigma_k,\delta_k$ to be recursively determined.\\
Now consider the actions on $\ue{B}_k$ and $\ue{B}_{k+1}$ of the linearized maps $\underline{\hat{\Psi}}$ and $\underline{\hat{f}}^{-n}$ respectively. Our purpose is to determine the windows centres and representative matrices in a way to satisfy the alignment test. We get
\begin{equation}\label{eq:w1w2quattrod}
\begin{array}{rclcrcl}
\ue{c}_1&=&\ue{X}_k'+A (\sigma_k,0,\delta_k,0)^T,&\qquad& \ue{c}_2&=&\ue{v}_n+G(n)\ue{p}_{k+1} \mbox{,}\\
W_1&=&A B_k,& & W_2 &=& G(n) B_{k+1} \mbox{.}
\end{array}
\end{equation}
Where $G(n)=\diag(L^{-n},1,L^n,1)$ and $\ue{v}_n=(0,-n \nu,0,0)$.\\
Recalling (\ref{eq:formageneralea}) then using (\ref{eq:formasempliceb}) and (\ref{eq:sigmaij}) the product $A B_k$ reads 
\begin{equation}\label{eq:prodottoab}
A B_k = \mu^{p}\left(
\begin{array}{cccc}
\mu a_{11} b_1 & \mu a_{12} b_2 & 0 & 0\\
0 &  b_2 & 0 & b_4\\
K_0(\mu) b_1 & \mu a_{32} b_2 & K_1(\mu)  b_3 & \mu K_2(\mu) b_4\\
\mu a_{41} b_1 & \mu a_{42} b_2 & 0 & 0
\end{array}
\right) \mbox{,}
\end{equation}
where
\begin{equation}
\begin{array}{rcl}
K_0(\mu) & := & -1/a+\mu a_{31}\\
K_1(\mu) & := \de &-\frac{a+\mu(a_{13}-a_{31}a^2)+\mu^2 a(a_{11}a_{33}-a _{13}a_{31})}{a(a+\mu a_{13})}\\
K_2(\mu) & := \de & \frac{a_{32}a-\mu(a_{12}a_{33}-a_{13}a_{32})}{a+\mu a_{13}}
\end{array}
\mbox{.}
\end{equation}

\subsubsection{Alignment test part one: intermediary matrices}
By using the obtained simplifications, the intermediary matrices read as
$$
M_k=\mu^{p}
\left(
\begin{array}{cccc}
\mu a_{11} b_1 &  \mu a_{12} b_2 & -c_{3}L^{-n} &0 \\
0 & b_2 & 0 & -c_4\\
K_0(\mu) b_1 & \mu a_{32} b_2 & -\mu L^n \Sigma_{33} c_3 & -\mu L^n \Sigma_{34} c_4\\
\mu a_{41} b_1 & \mu a_{42} b_2 & -\mu \Sigma_{43} c_{3} & -\mu \Sigma_{44} c_{4}\\
\end{array}
\right) \mbox{,}
$$
and 
$$
N_k=\mu^{p}
\left(
\begin{array}{cccc}
-c_1 L^{-n} & 0  & 0 & 0\\
0 & -c_2 & 0 & b_4\\
0 & 0 & K_1(\mu) b_3 & K_2(\mu) b_4\\
0 & 0 & 0 & 0
\end{array}
\right) \mbox{,}
$$
where 
\begin{equation}\label{eq:posizionibc}
b_i:=b_i^k,\qquad c_i:=b_i^{k+1} \mbox{,}
\end{equation}
so we have denoted by $c_i$ the entries of $B_{k+1}$.\\
As required by lemma \ref{lem:marco} we have to compute the inverse of $M_k$, however its explicit expression is quite complicate and cumbersome. We prefer a ``perturbative'' approach of this computation, by setting $L^n=\mu^{-\beta}$ and suitably choosing $\beta$. This leads to the following 
\begin{lem}
For all $p \in \mathbb{N}$, there exists $\beta_{p} \in \mathbb{N}$ such that, for all $\beta \geq \beta_{p}$, the matrix $M$ is invertible and its inverse admits the following expression
 \begin{equation}\label{eq:mmunoasintotica}
M_k^{-1}= \frac{\mu^{-p-1}}{2 D \Sigma_{33}}
\left(
\begin{array}{cccc}
 -\frac{\tilde{D}-a_{42} \Sigma_{33}}{b_1} & \mu \frac{a_{12} \tilde{D}}{b_1} & 0 & - \frac{a_{12} \Sigma_{33}}{b_1}\\
 -\frac{a_{41} \Sigma_{33}}{b_2} &  -\mu \frac{a_{11} \tilde{D}}{b_2} & 0 & \frac{a_{11} \Sigma_{33}}{ b_2}\\
  \frac{a_{41} \Sigma_{34}}{c_3} & \mu \frac{\Sigma_{34} D}{ c_3} & 0 &  \frac{a_{11} \Sigma_{34}}{c_3}\\
 -\frac{a_{41} \Sigma_{33}}{c_4} &  -\mu \frac{\Sigma_{33}D}{c_4} & 0 & \frac{a_{11} \Sigma_{33}}{c_4}
 \end{array}
\right)+O(\mu^{\beta-p-1}) \mbox{,}
\end{equation}
where $\tilde{D}:=\Sigma_{33}\Sigma_{44}-\Sigma_{34}\Sigma_{43}=D/\delta(\mu)$. Moreover, the explicit expression of the third column is
\begin{equation}\label{eq:terzacolonnamm1}
(M_k^{-1})_3=\frac{\mu^{\beta-p-1}}{2 D \Sigma_{33}}
\left(
\begin{array}{c}
-\frac{a_{12} \Sigma_{43}}{b_1}\\
\frac{a_{11}\Sigma_{43}}{b_2}\\
\frac{D-a_{11}\Sigma_{44}}{c_3}\\
\frac{a_{11} \Sigma_{43}}{c_4}
\end{array}
\right)+O(\mu^{2 \beta-p-2}) \mbox{.}
\end{equation}
\end{lem}
\nt As $\beta$ does not depend on $\mu$ this result holds at the price of a multiplicative $O(1)$ constant of the transition time $n(\beta):=\log_L \mu^{-\beta}=O(\log \mu^{-1})$.
\proof It is a straightforward check. By choosing 
\begin{equation}\label{eq:primalimitazbeta}
\beta_b:=4(p+1)
\end{equation}
one gets, for sufficiently small $\mu$ and $\beta\geq\beta_p$, 
\begin{equation}\label{eq:detmduattrod}
\det M_k  \sim -2 \mu^{4 p+3-\beta} b_1b_2 c_3 c_4 D \Sigma_{33}  \mbox{.}
\end{equation}
that is non-zero by (\ref{eq:splittingaij}). The adjugate can be computed in the same way, by using the dominance of terms containing $\mu^{\beta}$.
\endproof

\nt By using the already obtained form of $M_k^{-1}$ and by defining 
\begin{equation}
    \gamma:=\frac{a_{12}}{a_{11}} \mbox{,}
\end{equation}
the product $M_k^{-1}N_k$ yields, up to $O(\mu^{\beta})$
$$
M_k^{-1}N_k=\frac{1}{2}
\left(
\begin{array}{cccc}
0 & \gamma \frac{c_2}{b_1}  & 0 & -\gamma \frac{b_4}{b_1} \\
0 & -\frac{c_2}{b_2}  & 0 & \frac{b_4}{b_2} \\
0 & -\gamma \frac{c_2}{c_3}  & 0 & \gamma   \frac{b_4}{c_3} \\
0 &   \frac{c_2}{c_4}  & 0 & -\frac{b_4}{c_4} \\
\end{array}
\right) \mbox{.}
$$
Now, recalling (\ref{eq:posizionibc}) and reintroducing the dependencies on $k$, it is easy to see that if $b_{i}^k=b_{i}^{k+1}$ for some $k$, even in the simplest case $\gamma=0$, one gets
$$\norm{M_k^{-1}N_k}_{\infty}=\frac{1}{2}\left(1+\max \left\lbrace \frac{b_4^k}{b_2^k},\frac{b_2^k}{b_4^k}\right\rbrace \right) \geq 1 \mbox{,}$$
and the affine alignment test fails.\\
The situation is quite different in the anisochronous case, in which the $\theta$ component of the map $\ue{\hat{f}}^{-n}$ writes as (see \cite{mar96})
$$f_2^{-n}(\theta,\rho)=\theta-(\nu+\nu_1 \rho)n \mbox{.}$$
In this way, one gets an estimate of the form
$$\norm{M_k^{-1}N_k}_{\infty} \leq \frac{2}{|\alpha \nu_1 n-2|} \mbox{,}$$
(see \cite[Pag. 248]{mar96}), where $\alpha=O(\mu)$ and $n$ is kept sufficiently large.\\
Due to the presence of the term $\nu_1$, the quantity 
$\norm{M_k^{-1}N_k}_{\infty}$ tends to zero, as $n$ increases, and the alignment test can be satisfied.\\
This looks as the consequence of the ``transversality-torsion'' mechanism pointed out in \cite{cresguil1}: the joint action of splitting ($\alpha
\neq 0$) and torsion ($\nu_1 \neq 0$) transforms the partial hyperbolicity (i.e. in the two variables $(P,Q)$), in total hyperbolicity (in the four variables). In this way windows are compressed-stretched also in the $(\theta,\rho)$ variables.  This feature, already mentioned and used in \cite{east81}, is clearly stressed in lemma $4.3$ of the paper \cite{cres3}. In such result, eigenvalues of the map
$\ue{\mathcal{L}}:\mathcal{Q} \rw \mathcal{Q}$ defined by
$$\ue{\mathcal{L}}:=\ue{B}^{-1} \circ \ue{f}^{n} \circ \ue{\Phi} \circ \ue{B}$$
are explicitly computed, showing that for sufficiently small $\mu$ and large $n$, these are real and of modulus not equal to one. Consistently if $\alpha \nu_1 = 0$, the modulus of the eigenvalues corresponding to directions $(\theta,\rho)$ is equal to one and the phenomenon vanishes. \\
So the presence of torsion looks as an advantageous property in the WM machinery we cannot use. Nevertheless, by a careful choice of the representative matrices, as depending on $k$, it is possible to obtain a sort of ``simulation'' of the torsion effect:
\begin{lem}[simulated torsion]
Let $N=O(1/\mu)$. There exist an $O(1)$ positive constant $K$ and a sequence $ \{B_k \}_{k=1,\ldots,N} \in \mathbb{M}(4,4)$ of representative matrices of the form (\ref{eq:formasempliceb}) such that the resulting\footnote{in this computation the centres $\ue{p}_k$ does not appear and then these will be determined later.} $\ue{W}_{1,2}^a$ give rise to intermediary matrices satisfying  
\begin{equation}\label{eq:mmunon}
\norm{M_k^{-1}N_k}_{\infty} \leq 1-K \mu^2 \mbox{.}
\end{equation}
for all $k=1,\ldots,N$.
\end{lem}
As it is evident, this result gives a thin (but necessary) layer in order to satisfy the next estimates.
\proof
Suppose $\gamma \neq 0$, the simpler case $\gamma=0$ is an easy consequence.\\
Condition $\norm{M_k^{-1}N_k}_{\infty} < 1$ gives the following system of recursive inequalities
\begin{equation}\label{eq:sistemamnk}
\left\{
\begin{array}{rcl}
b_2^{k+1} + b_4^k & < & \frac{2}{|\gamma|} b_1^k\\
b_2^{k+1} + b_4^k & < & 2 b_2^k\\
b_2^{k+1} + b_4^k & < & \frac{2}{|\gamma|} b_3^{k+1}\\
b_2^{k+1} + b_4^k & < & 2 b_4^{k+1}\\
\end{array}
\right. \mbox{,}
\end{equation}
redefining $b_i^j:=|b_i^j|$, then the solution makes sense only if $b_i^k>0$ for all $k$.\\
Now note that as neither $b_1^{k+1}$ nor $b_{3}^k$ appear in the previous system, the choice of $b_{1}^k$ and of $b_{3}^k$ is free. Hence, the first and the third equations are trivially satisfied. In this way the previous system is satisfied (by a suitable choice of $b_1^k$ and $b_{3}^k$) for all $k$ if the following system 
\[
\left\{
\begin{array}{rcl}
b_2^{k+1}+ b_4^k & = & C b_2^k\\
b_2^{k+1}+ b_4^k & = & C b_4^{k+1}
\end{array}
\right.
\]
is, with $C<2$. By setting $x_k:=b_2^{k}$, the latter gives the following one dimensional discrete initial value problem
$$
\left\{
\begin{array}{rcl}
x_{k+2}&=&C x_{k+1}- x_k,\qquad k \geq 1 \\
x_2&=&b_2^2 \\
x_1&=&b_2^1
\end{array}
\right. \mbox{.}
$$
Fixed $N$, we are interested in a solution such that $x_k>0$ for all $k=1,\ldots,N$.\\
If $C<2$ the solution takes the form
\begin{equation}\label{eq:soluzricorrenza}
x_k=K_1 \cos (k-1) \alpha+ K_2 \sin (k-1) \alpha,\qquad k \geq 1 \mbox{,}
\end{equation}
where $\alpha=\arctan \frac{\sqrt{4-C^2}}{C}$.\\By choosing $x_2=(C/2)x_1$ we get $K_2=0$. Moreover, by setting $C:=\sqrt{4-\mu^2 b^2}$ where $b$ has to be determined, we have, $\alpha=\mu \frac{b}{2} + O(\mu^3)$. In this way (\ref{eq:soluzricorrenza}) reduces, up to higher orders, to 
\begin{equation}\label{eq:soluzricorr}
x_k=x_1 \cos \left(k \mu \frac{b}{2}\right) 
\end{equation}
for all $k=1,\ldots,N$. So $x_k>0$ for all $k=1,\ldots,N$ if $k \mu \frac{b}{2} \in (0,\pi/2)$, that is 
$$b=\frac{\pi}{\mu(N+1)} \mbox{.}$$
Now come back to the original set. Recalling the definition of the sequences $b_i^k$ in terms of $x_k$ we see that the remaining two inequalities of (\ref{eq:sistemamnk}) are satisfied by taking for all $k$
$$b_1^k=\frac{2}{C}\frac{\gamma}{2}(b_{2}^{k+1}+b_4^k)=\gamma  b_2^1,\qquad
b_3^{k+1}=b_1^k \mbox{,}$$ (we have used that $2/C>1$) for all $k$ ($b_3^1$ is arbitrary).\\
By keeping in mind that 
$x_1:=b_2^1$ is arbitrary, the remaining initial conditions are determined. The datum $b_4^1$ can be chosen as $b_4^1=(C/2) b_2^1$.\\
Note that $C=2-\frac{1}{4}\mu^2 b^2+O(\mu^4)$, so, by defining $K:=\frac{1}{4}b^2+O(\mu^2)$, we see that $b$ has to be an $O(1)$ constant, otherwise estimate (\ref{eq:mmunon}) and approximation (\ref{eq:soluzricorr}) do not properly make sense. This implies that $N=O(1/\mu)$, a fully compatible choice with our purposes.\\
\endproof
\begin{rem}
Note that variable parameters $b_{2,4}^k$ are involved in the $(\theta,\rho)$ part of the windows, while entries $b_{1,3}^k$, relative to the hyperbolic part can be chosen as constant during the evolution. In this sense we simulate the torsion effect.
\end{rem}
\begin{rem}
The decay of $b_{2,4}^k$ implies that the determinant of $B_k$, and then those of $M_k$,
approach to zero. Obviously, by hypothesis, zero is never attained but the entries $b_{2,4}^k$ reach values which are no longer $O(1)$. The smallest value of these entries (we get it for $k=N$) is immediate by (\ref{eq:soluzricorr})
\begin{equation}\label{eq:ultimiparametri}
b_2^N=b_2^1 \cos \left(\frac{\pi}{2}\frac{N}{N+1}\right) \sim
b_2^1 \frac{\pi}{2 N}=O(\mu) \mbox{.}
\end{equation}
In such a way (\ref{eq:detmduattrod}) changes as follows 
\begin{equation}\label{eq:detmmodificata}
\det M_k|_{k=N} \sim -2 \mu^{4 p+5-\beta} b_1^{N-1} b_2^{N-1} b_3^N
b_4^N D \Sigma_{33} \mbox{.}
\end{equation}
As an uniform estimates is handier, we extent condition (\ref{eq:primalimitazbeta}) to 
\begin{equation}\label{eq:limitazionebetaquattrodmod}
\beta_p := 6+4 p \mbox{,}
\end{equation}
for all $k$.
\end{rem}

\subsubsection{Alignment test part two: estimate of the term $M_k^{-1}(\ue{c}_2-\ue{c}_1)$}
Our aim is now to determine centres $\ue{p}_k$ in such a way the alignment test is satisfied for all $k$.\\
Set $\ue{\Delta}:=(\ue{c}_2-\ue{c}_1)$. By (\ref{eq:mmunon}), for a successful test we need, for all $k$
\begin{equation}\label{eq:mmunobma}
\norm{M_k^{-1}\ue{\Delta}}_{\infty} \leq \frac{1}{2} K \mu^2 \mbox{.}
\end{equation}
This means that centres of windows $\ue{W}_1^a$ and $\ue{W}_2^a$ should coincide ``as much as possible''. For this purpose we shall use the EEC mechanism in a very profitable way. 
\begin{lem}
Let us consider the sequence $B_k$ constructed before. For all ESC there exist an EEC, sequences $(\sigma_k,\delta_k)\in \R^2$ and $n_{k+1} \in \R$ in such a way the windows $\ue{W}_{1,2}^a$ (obtained by $\ue{B}_k$ and now completed by (\ref{eq:centriquattrod})) are correctly aligned for all $k$.\\
Furthermore, for all $k=1,\ldots,N-1$, the transition time from $\mathcal{T}_k$ to $\mathcal{T}_{k+1}$ satisfies 
\begin{equation}\label{eq:tempoditransizione}
n_{k+1}=O\left( \ln \frac{1}{\mu}\right)  \mbox{.}
\end{equation}
\end{lem}
\nt The previous formula, as $N=O(\mu^{-1})$, gives immediately the transition time estimate (\ref{eq:transitiontime}). 
\proof
We have to choose parameters $\delta_k$ and $\sigma_k$ in order to satisfy (\ref{eq:mmunobma}) for all $k$. It will be simpler to keep 
$\delta_{k+1}$ and $\sigma_{k+1}$ fixed and then determine $\delta_k,\sigma_k$.\\
Keeping in mind (\ref{eq:w1w2quattrod}), the two windows centres take the following form 
$$
\begin{array}{rcl}
\ue{c}_{1}&=&(\mu \sigma_k a_{11}+f_1(\mu) \delta_k, \theta_k', P_k'+\mu a_{33} \delta_k+K_0(\mu) \sigma_k,
\rho_{k+1}+\mu a_{41} \sigma_k + \mu a_{43} \delta_k)\\
\ue{c}_{2}&=&(\mu^{\beta}(Q_{k+1}+\sigma_{k+1}),\theta_{k+1}-n \nu, \mu^{-\beta}\delta_{k+1}, \rho_{k+1})
\end{array}
\mbox{.}
$$
In order to satisfy (\ref{eq:mmunobma}), we try to nullify as much components of $\ue{\Delta}$ as possible. This is achievable with first and fourth component via a suitable choice of $\delta_k$ and $\sigma_k$. More precisely, if
\begin{itemize}
\im $a_{41} \neq 0$, we choose $\sigma_k= -\delta_k a_{43}/a_{41}$ and then $\delta_k= \frac{\mu^{\beta}(Q_{k+1}+\sigma_{k+1})a_{41}}{f_1(\mu)a_{41}-\mu a_{11}a_{43}}$.
\im $a_{41}=0$, we take $\delta_k =0$ and then $\sigma_k= \mu^{\beta-1}(Q_{k+1}+\sigma_{k+1})/a_{11}$.
\end{itemize}
The described property holds for all $k=N-1,\ldots,1$. As free parameters, we choose 
$\delta_N=\sigma_N=0$. Note that $\delta_k$ can be either zero or O($\mu^{\beta})$. This implies that $$\Delta_3=O(1) \mbox{.}$$Moreover, if $\beta$ is sufficiently large, the third column of $M_k^{-1}$ is arbitrarily small and the contribution of the product $M_k^{-1}\Delta_3$ is not greater than $O(\mu^{3(p+1)})$ and then negligible for suitable $p$.\\
Now we are going to use elasticity of the transition chain in order to nullify $\Delta_2$. For this purpose, recall the notational setting of sec. \ref{subsec:sequences} and consider values $y_{k+1}$ and $y_k$ of a given ESC. By keeping $\tilde{y}_{k+1}$ fixed we want to move $y_k$ to a suitable $\tilde{y}_k \in E_k$.\\
Note that for all $y \in E_k$, the argument we have used to get $\Delta_1=\Delta_4=0$ can be repeated simply by replacing the coordinates of $\ue{X}_k,\ue{X}_k'$ with those of $\ue{X}_k(y),\ue{X}_k'(y)$. This is possible as the elements of the family $k+1$ do not change and the algorithm remains well posed. Now $\delta_k,\sigma_k$ looks as functions $\delta_k(y),\sigma_k(y)$.\\
Define $2e_k:=\diam E_k$ and
$$
\begin{array}{rcl}
\ml{O}_{k}^-&:=&\{\theta_{k}'(y): \quad y \in
(y_{k}-e_{k},y_{k}]\}\\
\ml{O}_{k}^+&:=&\{\theta_{k}'(y): \quad y \in
[y_{k},y_{k}+e_{k})\}
\end{array}
\mbox{.}
$$
In the time $n(\beta)$ the section $\Sigma_N$ is visited by points whose reciprocal distance is not greater than  $O((n(\beta))^{-\frac{1}{\tau}})$ (see e.g. \cite{gal99} and references therein). As $\diam \ml{O}_{k}^{\pm}=O(1)$ there exist subsequences $n_{\pm}^{m}$ with $n_{\pm}^1$ of the same order of $n(\beta)$, in a way that 
$$\theta_{k+1}-n_{\pm}^m \nu \in \ml{O}_{k}^{\pm} \mbox{.}$$
Suppose that $n_{k+1}:=n_+^{\bar{m}}$ for some chosen $\bar{m}$ and define
\begin{equation}
	\ml{F}(y):=[\theta_k'(y)-(\theta_{k+1}-n_{k+1}\nu)]\modulus 2 \pi \equiv \Delta_2 \mbox{.}
\end{equation}
From (\ref{eq:variabilitaphiuno}) we deduce that $\ml{F}(y)$ is continuous and $\frac{d}{dy} \ml{F}(y) =O(\mu^{-1})$ on $E_k$, moreover $\ml{F}(y_k) \ml{F}(y_k+e_k) <0$, so there exists unique\footnote{note that this point can be constructively determined with an arbitrary precision.} $\tilde{y}_k \in E_k$ such that $\ml{F}(y_k)=0$.\\
In such a way $\norm{M_k^{-1}\ue{\Delta}}_{\infty} \sim 0$ and (\ref{eq:mmunobma}) is trivially satisfied. 
\endproof

\begin{rem}
The elasticity of the transition chain could seems as an unnecessary tool to achieve affine alignment. Without using it we could obtain 
\begin{equation}
    \norm{M_k^{-1}\ue{\Delta}}_{\infty} \leq \mu^{-p} J_k  \Delta_2 \mbox{,}
\end{equation}
where $J_k= \max\{\frac{\gamma}{b_1^k},\frac{1}{b_2^k},\frac{\gamma}{b_{3}^{k+1}},\frac{1}{b_{4}^{k+1}}\}$. Clearly, due to the decreasing of the sequences $b_{2,4}^k$, the greatest value of $J_k$ is $J_N=O(\mu^{-1})$. In order to satisfy (\ref{eq:mmunobma}), $\Delta_2$ should be of an order not greater than $\mu^{3+p}$, requiring in this way the following ergodization time 
$$ n =O (\mu^{-(p+3)}) \mbox{.}$$
As it will be clear by looking at the estimates of the next section, during such (long) time the remainder of $\underline{f}^{-n}$ is not suitably bounded and the affine problem is not still a good approximation of the complete one.
\end{rem}
\begin{rem}
The use of the EEC allows us to bypass the ergodization time on $\Sigma_N$ so the drift time is independent on the Diophantine constants, as in \cite{berbol}.
\end{rem}

\subsection{Estimate of the remainders}
Recall lemma \ref{lem:marco}. By construction, functions $\hat{W}_1(\ue{x})$ and $\hat{W}_2(\ue{x})$ are given by remainders (\ref{eq:restof}) and (\ref{eq:restopsi}) respectively, evaluated at $\ue{\xi} \leftarrow (\sigma_k,0,\delta_k,0)+B_k \ue{x}$, we denote with $\xi(\ue{x})$, $\ue{x} \in (-2,2)^4$.\\ 
Furthermore, set $\gamma_k:=\norm{M_k^{-1}}_{\infty}$ and define 
\begin{equation}\label{eq:chiunochidue}
\begin{array}{rcl}
\chi_1 & := & \sup_{\underline{x} \in \mathcal{U}}[\norm{M_k^{-1} \underline{R}_{k+1}(\underline{x})}_{\infty}+\gamma_0 \norm{\underline{\hat{R}}_k(\underline{x})}_{\infty}]\\
\chi_2 & := & \sup_{\underline{x} \in \mathcal{U}}[\norm{M_k^{-1} \pl_{\underline{x}} \underline{R}_{k+1}(\underline{x})}_{\infty}+\gamma_0
\norm{\pl_{\underline{x}}\underline{\hat{R}}_k(\underline{x})}_{\infty}]
\mbox{.}
\end{array}
\end{equation}
In this way
$$\chi^c \leq \max\{\chi_1,\chi_2\} \mbox{.}$$
First of all, it is easy to see from (\ref{eq:restopsi}) that
$$
\norm{\hat{\ue{R}}_k(\ue{x})}_{\infty} \leq C_2 \norm{\xi(\ue{x})}_{\infty}^2,\qquad \norm{\pl_{\ue{x}}\hat{\ue{R}}_k(\ue{x})}_{\infty} \leq C_3 \norm{\xi(\ue{x})}_{\infty} \norm{\pl_{\ue{x}}\xi(\ue{x})}_{\infty} \mbox{,}
$$
where $C_{2,3} < +\infty$ by regular (analytic) dependence on initial data over finite time. Moreover, by a comparison between (\ref{eq:mmunoasintotica}) and (\ref{eq:ultimiparametri}), we have $\max_k \gamma_k=O(\mu^{-1-p})$. So we get
$$
\max_k \gamma_k \norm{\hat{\ue{R}}_k(\ue{x})}_{\infty} \sim 
\max_k \gamma_k \norm{\pl_{\ue{x}}\hat{\ue{R}}_k(\ue{x})}_{\infty}=O(\mu^{p-2}) \mbox{,}
$$
uniformly in $\ue{x}$.\\
By looking at formula (\ref{eq:restof}), we obtain 
$$
\begin{array}{rcl}
\max_k \norm{M_k^{-1}\ue{R}(\ue{x})}_{\infty} &\sim &
\mu^{-1-p} \left\| \left(
\begin{array}{cccc}
m_{11} & \cdots & m_{13} \mu^{\beta-1} & \cdots\\
m_{21}\mu^{-1} & \cdots & m_{23} \mu^{\beta-2} & \cdots\\
m_{31} & \cdots & m_{33} \mu^{\beta-1} & \cdots\\
m_{41}\mu^{-1} & \cdots & m_{43} \mu^{\beta-2} & \cdots
\end{array}
\right)
\left(
\begin{array}{c}
n \mu^{\beta-p}\\
0\\
n \mu^{- \beta+2 p}\\
0
\end{array}
\right) \right\|_{\infty} \\
&=& O(n_{k+1} \mu^{p-3})
\end{array}
$$
recalling (\ref{eq:limitazionebetaquattrodmod}) and (\ref{eq:tempoditransizione}). Similarly 
$$
\max_k \norm{M_k^{-1} \pl_{\ue{x}} \ue{R}_k(\ue{x})}_{\infty}=O(n_{k+1}
\mu^{p-3}) \mbox{.}
$$
Keeping in mind (\ref{eq:mmunon}), in order to satisfy (\ref{eq:condizionechic}b) (and then (\ref{eq:condizionechic}a)), it is sufficient to choose $p=6$ (from which $\beta \geq 30$). As $\lim_{\mu \rw 0} n_{k+1} \mu =0$ by definition, we have $\chi_{1,2}=o(\mu^2)$ and the proof is complete.
 
\subsection*{Acknowledgements} This work is part of \cite{for}, I wish to thank my advisor Prof. G. Gallavotti and Prof. G. Gentile for all the helpful discussions, advices and precious comments on it.
I am grateful to all the specialists for having kindly answered to my questions on their papers with stimulating comments, and in particular to Proff. M. Berti, L. Biasco, L. Chierchia, J. Cresson and J.P. Marco.

\bibliographystyle{alpha}
\bibliography{Fdad}

\end{document}